\input amstex
\input epsf
\documentstyle{amsppt}

\def\rmop#1{\expandafter\def\csname#1\endcsname{\operatorname{#1}}} \rmop{diam}
\rmop{sh} \rmop{lsh} \rmop{dim} \rmop{id} \rmop{N} \rmop{Int}
\let\but\setminus \def\Cl#1{\overline{#1}} \def\R{\Bbb R} \def\heps{{\eps/2}}
\TagsOnRight \let\eps\varepsilon \let\dta\delta \let\phi\varphi \let\tl\tilde
\let\emb\hookrightarrow \let\su\subset \let\x\times \redefine\o{\circ}
\redefine\sp#1{@, @, @, @, @, \langle @! @! @! @! @! #1\rangle}
\def\col{{\ssize~\searrow~}} \def\suchthat{{\bigm|\,}}

\topmatter 

\title Singular link concordance implies link homotopy \\
in codimension $\ge 3$
\endtitle
\rightheadtext{Singular link concordance implies link homotopy}
\author S. Melikhov \endauthor

\abstract\nofrills
We prove the analogue of the Concordance Implies Isotopy in Codimension
$\ge 3$ Theorem for link maps, together with some other its singular analogues.
\endabstract



\endtopmatter

\document 

\head 1. Introduction and statement of results \endhead

An interesting and deep question of geometric topology is to determine under
which conditions two given maps from a special class (for example link maps,
immersions, embeddings, close embeddings) are connected by a continuous
deformation in this class (respectively by link homotopy, regular homotopy,
isotopy, small isotopy).
Sometimes this problem of map {\it equivalence} can be reduced to a problem of
{\it existence} of (extensions of) maps of the same class.
A reduction of isotopy to relative embeddability was achieved in 60s by
Smale and Haefliger \cite{10, 1.2}, Zeeman \cite{38}, Lickorish \cite{23,
Th\. ~6} and Hudson \cite{12}, \cite{13}, who showed that concordant embeddings
are ambient isotopic in codimension $\ge 3$, in smooth and PL categories
(see \cite{30}, \cite{27}, \cite{2, \S7} for alternative proofs, and
\cite{31, last historic remark}, \cite{11, proof of Cor\. ~1} for typical
applications of this reduction).
The similar question for link maps (does singular link concordance imply link
homotopy in codimension $\ge 3$?) has remained open \cite{21, p\. ~303}.

Throughout this paper, let $X=X_1\sqcup\dots\sqcup X_k$ be a compact polyhedron
($X_i$'s are fixed, but not necessarily connected) and $Q$ a PL manifold.
A continuous map $f\:X\to Q$ is called a {\it (generalized) link map} if
$fX_i\cap fX_j=\emptyset$ whenever $i\neq j$.
A {\it (singular) link concordance} -- not to be confused with (embedded)
concordance -- between link maps $f_0,f_1\:X\to Q$ is a link map
$F\:X\x I\to Q\x I$ such that $F(x,i)=(f_ix,i)$ for $i=0,1$ and each $x\in X$.
A {\it level-preserving} (i.e\. such that $F(X\x t)\su Q\x t$ for each
$t\in I$) link concordance is called a {\it link homotopy}.

Under link concordance and link homotopy we also mean the equivalence relations
they generate on the set of link maps $X\to Q$.
In this section all maps are supposed to be continuous unless the contrary is
stated. Isotopy means ambient isotopy.

The relation of link homotopy was introduced by Milnor for classical links of
circles in $\R^3$ \cite{26}.
Scott considered spherical link maps (i.e\. of spheres into a sphere) in higher
dimensions up to link homotopy \cite{33}, and Nezhinskij -- up to link
concordance (see \cite{28, \S1}).
Since mid-80s the problem of classification of spherical and generalized
link maps up to link homotopy and link concordance has been studied widely
(see for instance \cite{4}, \cite{9}, \cite{16}--\cite{22}, \cite{28},
\cite{29}, \cite{32}, \cite{35}).

One may have an intuitive feeling for link maps up to link homotopy as for
`links modulo knots', that is, ignoring knotting phenomena and concentrating
on interaction of different components, but it is well also to bear in mind
the following.
Some embedded links with unknotted components, which are link homotopic to the
trivial link, are not isotopic to it (e.g\. the Whitehead link or the link on
Fig\. ~1a; see \cite{25, 7.7}, \cite{7, \S2} for higher dimensional examples).
Also, leaving alone the case of link maps $X^n\to Q^m$ with negative
codimension $m-n$, observe that there are link maps not link homotopic to
embedded links (see Fig\. ~1b and codimension-two examples in \cite{4},
\cite{16}, \cite{17, 2.22}, \cite{22, \S3,4}, \cite{25, 7.9}).

\bigskip
\def\epsfsize#1#2{0.7\hsize}
\centerline{\epsffile{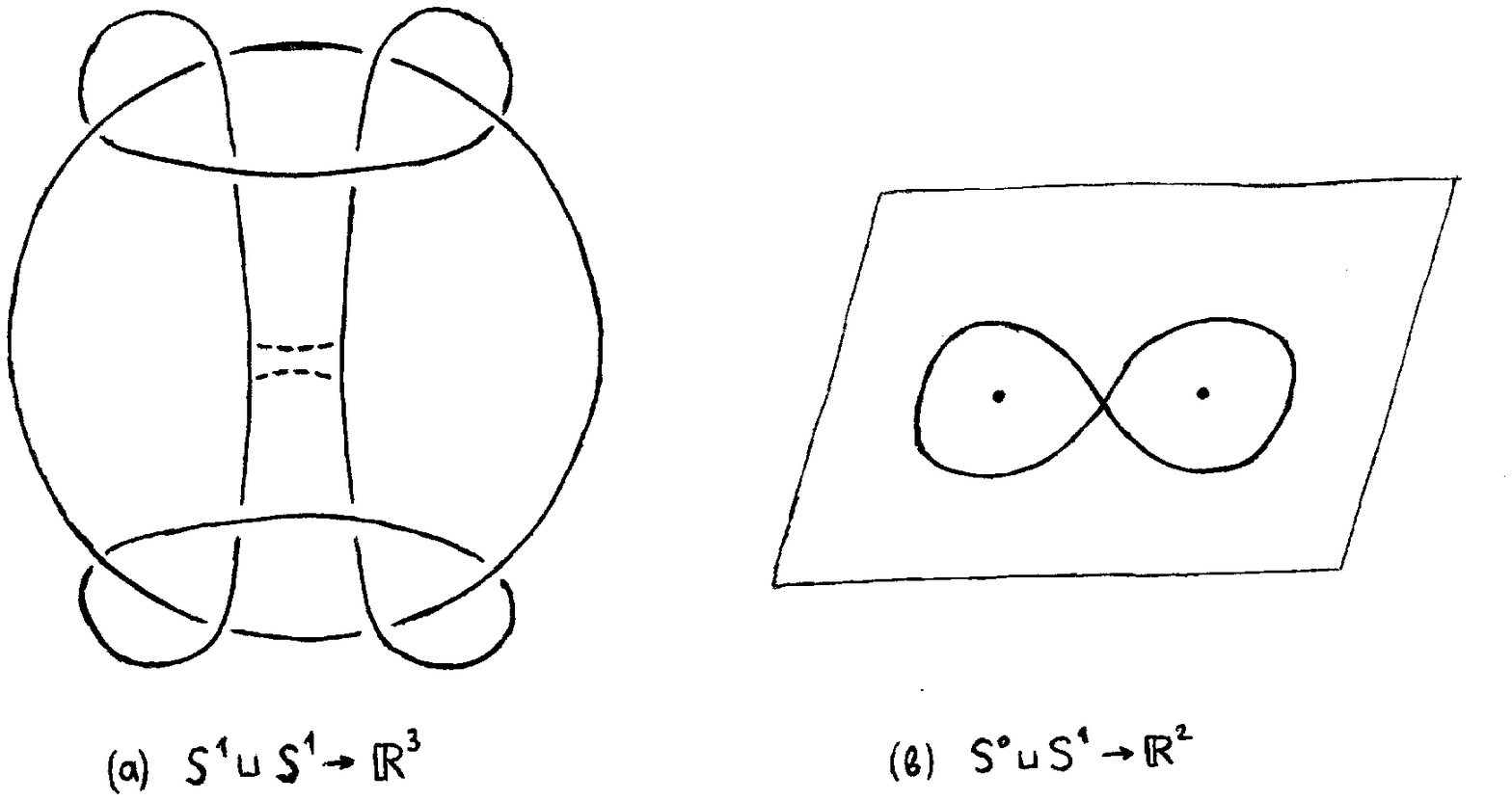}}
\medskip
\centerline{\bf Figure 1}
\bigskip

It appears that for classical links $S^1\sqcup\dots\sqcup S^1\emb S^3$
in PL category embedded concordance implies link homotopy \cite{5}, \cite{6}.
Notice that it does not imply isotopy, see Fig\. 1a: to obtain an unlinking
concordance, make a bridge as indicated by dotted line, isotop one of two
appeared circles far from others and glue it up by a disk, cf. \cite{37}.
Moreover, singular link concordance implies link homotopy for link maps
$S^1\sqcup\dots\sqcup S^1\to S^3$ \cite{24} (see also \cite{8}, \cite{36}).

Does link concordance imply link homotopy in general? Koschorke proved that it
does for link maps $S^{n_1}\sqcup\dots\sqcup S^{n_k}\to S^m$, where
$n_1\le m-3$ and $n_i\le\frac{2}3m-1$ for $i>1$ \cite{20}, and
for smooth embedded links $S^{n_1}\sqcup\dots\sqcup S^{n_k}\emb S^m$, where
$n_i\le m-3$ for all $i$, $n_i\le l(m-2) - (n_1+\dots+n_k)$ for $i>1$, and
all strict sublinks are link homotopically trivial \cite{21, 8.4c}.
For two-component link maps $S^p\sqcup S^q\to S^m$ link concordance was shown
to coincide with link homotopy in the range $2p+2q\le 3m-5$ \cite{9}.
Teichner announced that link concordance implies link homotopy for spherical
link maps in codimension $\ge 2$ \cite{36, Rem\. to Lem\. 3.1},
\cite{21, footnote on p\. ~303}.
Meanwhile Sayakhova obtained an example of link concordant but not
link homotopic link maps $S^1\sqcup S^1\sqcup S^2\to S^3$ \cite{32}.
In this paper we prove the following

\proclaim{Theorem 1.1} Let $X^n=X_1\sqcup\dots\sqcup X_k$ be a compact
polyhedron, $Q^m$ a PL manifold, $m-n\ge 3$, and $f_0,f_1\:X\to Q$ link maps.
If $f_0$ and $f_1$ are link concordant, then they are link homotopic.
\endproclaim

Thus in codimension $\ge 3$ any classification of link maps up to link
concordance (in particular, ones in \cite{29}, \cite{35}) is that up to link
homotopy, and any invariant of link homotopy is that of link concordance.
In contrast, the sufficiency of link concordance for link homotopy in \cite{9},
\cite{20}, \cite{21} was obtained as a corollary of completeness and link
concordance invariance of certain link homotopy invariants.

To prove Theorem 1.1 we use a different approach, independent on the recent
progress in the theory of link maps.
The proof goes in PL category, and the tool we use for codimension $\ge 3$
is a version of Zeeman's `sunny collapsing'.
Actually, our method develops some ideas from a part \cite{13, 5.1} (see also
\cite{12, 9.5}) of Hudson's proof of the Concordance Implies Isotopy Theorem.
However, the proof of 1.1 seems to be the singular analogue of a new proof
(to appear in a subsequent paper) of the Concordance Implies Isotopy Theorem
rather than of Hudson's or of either alternative proof mentioned above.
Apparently, the ideas of this paper do not suffice to prove 1.1 in codimension
2, where Concordance Implies Isotopy fails.
To make a brief introduction into the proof of 1.1, we sketch its idea later
in this section.

Given link maps $f_0,f_1\:S^{n_1}\sqcup\dots\sqcup S^{n_k}\to S^m$, where
$n_i\le m-2$, one defines a (set-valued) {\it connected sum}
$f_0\sharp f_1\:S^{n_1}\sqcup\dots\sqcup S^{n_k}\to S^m$ as follows.
Choose a basepoint $a_i$ in each $S^{n_i}$.
Push the image of $f_0$ into the southern hemispere of $S^m$, so that each
component $f_0(S^{n_i})$ meets the equator precisely in one point $f_0(a_i)$,
and push similarly the image of $f_1$ into the northern hemisphere so that,
in addition, $f_1$ agrees with $f_0$ on $a_i$'s (this is easy to achieve
using general position).
Now define $f_0\sharp f_1$ by shrinking the equator of each $S^{n_i}$ onto
$f_0(a_i)=f_1(a_i)$ and mapping the southern and northern hemispheres using
(shifted) $f_0$ and $f_1$ respectively.

Connected sum induces an operation (also called connected sum) on the set
$LM^m_{n_1,\dots,n_k}$ of link homotopy classes of link maps
$S^{n_1}\sqcup\dots\sqcup S^{n_k}\to S^m$ (by picking up arbitrary
representatives of link homotopy classes).
This operation is known to be well-defined and single-valued when
$n_i\le m-3$ \cite{19}, \cite{25}, \cite{33}.
Generally, this is not the case in codimension two \cite{22, Fig\. ~4.1},
except for some special situations, such as two-component \cite{17, 2.3}
and base-point preserving \cite{19, 1.4} link maps (although in the latter
case connected sum fails to be commutative \cite{19, 3.12}).

A {\it reflection} $(-f)$ of a link map $f\:S^{n_1}\sqcup\dots\sqcup S^{n_k}\to
S^m$ is the composition $R\o f\o r$, where $r$ and $R$ are reflections of
$S^{n_i}$'s and of $S^m$ in their equators.
Clearly, reflections of link homotopic link maps are link homotopic (by the
`reflected' link homotopy).
Thus, reflection induces a well-defined operation (also called reflection)
on $LM^m_{n_1,\dots,n_k}$.

\proclaim{Corollary 1.2} Suppose that $n_i\leq m-3$ for each $i=1,\dots,k$.
Then connected sum (regarded as addition) and reflection (regarded as
inverse) generate abelian group structure on $LM^m_{n_1,\dots,n_k}$.
\endproclaim

In fact, a statement similar to 1.2 can be found in \cite{33}.
But, as it was pointed out in \cite{17, 2.4}, the proof of this statement
had contained a gap: it had been actually shown that reflection is inverse up
to link concordance, not up to link homotopy.
Theorem 1.1 covers this gap; the rest of the proof of 1.2 goes as in
\cite{33}, \cite{19, 1.4} (see also \cite{10}, \cite{17, \S2},
\cite{18, 1.19}, \cite{23, Th\. ~7}, \cite{25, 5.2}, \cite{28}).

Corollary 1.2 generalizes \cite{17, 2.3}, \cite{18, 1.19}, \cite{19, 1.4+1.7}
and answers \cite{9, Question ~1} in the case of codimension $\geq 3$.
Compare 1.2 to group structures on link maps $S^p\sqcup S^q\to S^m$ up to link
concordance \cite{28}, on embedded links up to link homotopy \cite{25}, on
smooth embeddings of spheres up to diffeotopy \cite{10}.

Analogues of link maps, called `doodles', were introduced by Fenn and Taylor
in `search for a method of cancelling triple points' \cite{3} (see also
\cite{15}).
We call a continuous map $f\:X_1\sqcup\dots\sqcup X_l\to Q$ a
{\it (generalized) doodle} if $fX_i\cap fX_j\cap fX_k=\emptyset$ for any
distinct $i,j,k$.
In the obvious way one can define {\it doodle concordance} and {\it doodle
homotopy} (or, in terms of \cite{3}, cobordism and isotopy of doodles).
See \cite{3} for examples of doodles $S^1\sqcup\dots\sqcup S^1\to\R^2$,
in particular, of a doodle which is doodle concordant but not doodle homotopic
to the trivial doodle. In codimension $\ge 3$, our method works to
prove the following

\proclaim{Theorem 1.3} Let $X^n=X_1\sqcup\dots\sqcup X_k$ be a compact
polyhedron, $Q^m$ a PL manifold, and $f_0,f_1\:X\to Q$ doodles.
If $f_0$ and $f_1$ are doodle concordant, then they are doodle homotopic,
provided $m-n\ge 3$.
\endproclaim

One can define {\it $l$-doodle} to be a map $f\:X_1\sqcup\dots\sqcup X_k\to Q$
such that images of any $l$ distinct $X_i$'s have no point in common. Then
`$l$-doodle concordance implies $l$-doodle homotopy' theorem is stated and
proved analogously.
We prove the following more general statement on making a map preserve
levels almost without introducing new singularities:

\proclaim{Theorem 1.4} Let $X^n$ be a compact polyhedron, $Q^m$ a PL manifold,
$m-n\ge 3$, and $F\:X\x I\to Q\x I$ a map such that $F(X\x i)\su Q\x i$,
$i=0,1$.
For each $\eps>0$ there is a level-preserving map $\Phi\:X\x I\to Q\x I$ such
that $\Phi(x,i)=F(x,i)$ for $i=0,1$ and each $x\in X$, and such that for any
$x_1,\dots,x_l\in X$ the following holds:
$\Phi(x_1\x I)\cap\dots\cap\Phi(x_l\x I)\neq\emptyset \text{\ \ only if\ \ }
F(\N_\eps(x_1)\x I)\cap\dots\cap F(\N_\eps(x_l)\x I)\neq\emptyset.$
\endproclaim

Here $\N_\eps(x)$ denotes the $\eps$-neighborhood of $x$ in some fixed metric
on $X$.
Theorem 1.4 is proved in essentially the same ideas as its partial cases 1.2
and 1.3.
Notice that the condition on $\Phi$ in 1.4 is equivalent to the requirement
that, for each positive integer $l$, the `singular set' of $\Phi$ lies in an
arbitrarily small neighborhood of the `singular set' of $F$, both `singular
sets' being considered as subsets of the configuration space $X\x\dots\x X$ of
$l$-tuples of generators $x\x I$.
The `singular set' of $F$ always contains the diagonal of $X\x\dots\x X$,
hence even if $F$ is an embedding, $\Phi$ may be not an embedding.
Therefore the (PL/smooth) Concordance Implies Isotopy Theorem is not a
partial case of (the PL/smooth version of) Theorem 1.4. Instead, we have
Corollary 1.5 below, which can be thought of as the singular version of
Concordance Implies Isotopy.

It is said that $X$ is {\it quasi-embeddable} in $Q$, if for each $\eps>0$
there exists an {\it $\eps$-map} $f_{\eps}\:X\to Q$, i.e\. such map that
$\diam f_{\eps}^{-1}(q)<\eps$ for each $q\in f(X)$.
Quasi-embeddability and embeddability of $X^n$ into $Q^m$ are equivalent in
case $m\ge\frac{3(n+1)}2$ by \cite{11, Th\. 1}, and distinct for all pairs
$(m,n)$ such that $3<m<\frac{3(n+1)}2$, see \cite{34}.
Call embeddings $f_0,f_1\:X\emb Q$ {\it quasi-concordant (quasi-isotopic)} if
for each $\eps>0$ there is a (level-preserving) $\eps$-map $F\:X\x I\to Q\x I$
with $F(x,i)=(f_ix,i)$.
Clearly, quasi-isotopy is stronger than link homotopy, but weaker than isotopy
(consider the link on Fig\. ~1a and the trivial link).
If $m>\frac{3(n+1)}2$, quasi-isotopy implies PL isotopy for PL
embeddings $X^n\emb Q^m$ by \cite{11, Cor\. ~1}.

\proclaim{Corollary 1.5} Let $X^n$ be a compact polyhedron, $Q^m$ a PL
manifold, $m-n\ge 3$, and $f_0,f_1\:X\emb Q$ embeddings.
If $f_0$ and $f_1$ are quasi-concordant, then they are quasi-isotopic.
\endproclaim

Any map $f\:X\to Q$ is homotopic, by an arbitrarily small homotopy, to a PL map
and, moreover, the homotopy can be chosen to fix any subpolyhedron $Y$ such
that $f|_Y$ is PL \cite{39}.
Therefore link concordant link maps $f_0,f_1\:X\to Q$ are link homotopic to PL
link concordant PL link maps $f_0',f_1'\:X\to Q$.
Thus Theorem 1.1 is reduced to its PL version (analogously for Theorems 1.3,
1.4), and what we prove below is these PL versions.
Actually, 1.1 -- 1.5 are equivalent to their PL versions (for level-preserving
approximations can be made) and to their smooth versions (for similar
approximations in smooth category can be made).

Let us outline the idea of proof of Theorem 1.1.
Let $F\:X\x I\to Q\x I$ be the given link concordance between $f_0$ and $f_1$.
Suppose that sun emits its rays along the $I$-fibers of $Q\x I$, upside down.
The idea is to (singularly) reparametrize the product $X\x I$ by means of a map
$H\:X\x I\to X\x I$ (fixing $X\x \{0,1\}$, but not necessarily fiber-preserving
or homeomorphism), so that the $F$-image of each `fiber' $H(X\x t)$ is not
self-overshadowing (although possibly self-intersecting).
Then simple vertical shifting of each $F\o H(X\x t)$ into $Q\x t$ introduces
no new singularities, while the levels get preserved under shifted $F\o H$
(compare this to a proof of link concordance invariance of the
$\alpha$-invariant in \cite{35}).

To obtain such reparametrization, notice first that without loss of generality
$F$ is PL and by general position (\S2) we can also assume that the set $S$ of
points in $X\x I$, which $F$-images lie in the same sunray with some other ones,
is of codimension $\ge 2$ in $X\x I$.
Hence $Y\x I\but S$ is connected for each connected component $Y$ of $X$, and
it follows (\S3) that we can collapse $X\x I$ onto $X\x 0$ so that the
simplexes of $S$ are collapsed in any order of decreasing dimension, in
particular, in the order their $F$-images overshadow each other.
In other words, there is a collapse $X\x I\col X\x 0$ with the following
`sunny' property: `$F(a)$ overshadows $F(b)$' implies `$a$ is collapsed before
$b$' for any points $a,b\in X\x I$.
By a special care (\S4; see remark in the next paragraph) this collapse can
be improved to satisfy the following `stable sunny' condition: `$F(a)$
overshadows $F(b)$' implies `$a$ is collapsed before a neighborhood of $b$'
for any points $a,b\in X\x I$.
This slight improvement is, however, the key step.
Indeed, by it we achieve that the $F$-image of the frontier (in $X\x I$)
of points, already collapsed at the moment, is not self-overshadowing (for
each moment $t\in I$ during the collapse).
Thus the mapping of $X\x I$ onto itself, defined by $(x,t)\mapsto$
`the $t$-moment image of $x\x 1$ under the collapse' (which lies in the
$t$-moment frontier), is a proper reparametrization.

Actually, only (allowable) new self-intersections of each component
$F\o H(X_i\x t)$ can appear, if we weaken the conditions of reparametrization
and of stable sunny collapse so that shadows, casted by each $F(X_i\x I)$ onto
itself, are not taken into account.
To simplify the proof, we obtain only these weakened properties. The same
method works to prove 1.3, while in the proof of 1.5 we do not take into
account the shadows between sufficiently close points only.

Hereafter we assume all spaces to be polyhedra and all maps to be
piecewise-linear, unless the contrary is stated.
We follow \cite{31} for the PL notation.
We denote simplicial complexes and their bodies by the same letters.

\head 2. General position \endhead

\definition{Definition} We call $S(f)=\Cl{\{a\in A\suchthat
f^{-1}f(a)\neq a\}}$ the {\it singular set} of a PL map $f\:A\to B$
(not to be confused with the `singular set' in the configuration space,
mentioned in \S1; hereafter we use only the last definition).
A subpolyhedron $L$ of a polyhedron $K$ is {\it locally of codimension $\ge k$}
in $K$, if every $n$-simplex of $L$ faces some $(n+k)$-simplex of $K$ for any
triangulation of the pair $(K,L)$ (cf\. \cite{23}).
Let us say that a map $F\:X\x I\to Q\x I$ is a {\it general position map} with
respect to the projections $P\:Q\x I\to Q$, $p\:X\x I\to X$ if

1) $S(P\o F)$ is locally of codimension $\ge m-n-1$ in $X\x I$, and

2) $p|_{S(P\o F)}$ is non-degenerate.
\enddefinition

\proclaim{Lemma 2.1}
Let $X^n$ be a compact polyhedron, $Q^m$ a PL manifold, $m-n\geq 2$ and
$F\:X\x I\to Q\x I$ a PL map.
For each $\eps>0$, $F$ is $\eps$-homotopic to a PL general position map
$\tl F\:X\x I\to Q\x I$ such that $\tl F^{-1}(Q\x t)=F^{-1}(Q\x t)$ for each
$t\in I$.
\endproclaim

\demo{Proof}
For each point $p\in X\x I$, we fix the $I$-coordinate of $F(p)$ and change
its $Q$-coordinate, following Bing \cite{1, proof of 2.1} and assuming
the following remarks.
The property (2), although not included into the statement of \cite{1, 2.1},
was actually achieved in its proof (see \cite{1, 2.5}).
The additional restriction $m-n\geq 3$ in \cite{1, 2.1} was used only to
prove some properties additional to (1) and (2), hence it can be omitted here.
(However, in the sequel we use only the codimension $\ge 3$ case of Lemma 2.1.)
\qed \enddemo

The analogous to 2.1 result for embedding $F$, with condition (2) dropped,
was obtained in \cite{14, Lem\. ~1} (compare to \cite{13, 5.2}).

\head 3. Sunny collapsing \endhead

\definition{Definition} Let us think of the second factor of $Q\x I$ as of
heigth (that is, a point $(q_1,t_1)$ lies {\it below} a point $(q_2,t_2)$ if
$q_1=q_2$ and $t_1<t_2$). If $V\su Q\x I$, let $\sh V$ denote {\it shadow}
of $V$, the set of points of $Q\x I$ lying (strictly) below some point of $V$.
We say that a collapse $V\col W$ in $Q\x I$ is a {\it simple sunny} collapse,
if no point of $V\but W$ lies in $\sh V$ (equivalently, $V\cap\sh V\su W$).
A sequence of simple sunny collapses is called a {\it sunny} collapse (this
concept is due to Zeeman \cite{38}). It is convinient to generalize this
for the case of any map $F\:K\to Q\x I$, $K$ being a polyhedron. For any
$V\su K$ define {\it $F$-shadow} of $V$ by $\sh_F V=F^{-1}(\sh F(V))$.
A collapse $V\col W$ in $K$ is called a {\it simple $F$-sunny} collapse if
$V\cap\sh_F\su W$, and an {\it $F$-sunny} collapse if it is a sequence of
simple $F$-sunny ones.
\enddefinition

\proclaim{Lemma 3.1} Let $X^n$ be a compact polyhedron, $Q^m$ a PL manifold and
$m-n\geq 3$. If $F\:X\x I\to Q\x I$ is a PL general position map such that
$F(X\x 0)\su Q\x 0$, then there is an $F$-sunny collapse $X\x I\col X\x 0$.
\endproclaim

For the case of embedding $F$ this was proved in \cite{14, Lem\. ~2}
(the condition $F(X\x 0)\su Q\x 0$ was erroneously omitted in the statement).
The proof goes in general case with only minor changes; nevertheless for
completeness we sketch it below (see also \cite{38, proof of Lem\. ~9}
for detailed proof of a similar statement).

\demo{Proof} Assume that Lemma 3.1 is true when $\dim X<n$ and prove it for
$\dim X=n$ (if $n=0$ then $S(P\o F)=\emptyset$ and any collapse
$\{x_1,...,x_m\}\x I\col\{x_1,...x_m\}\x 0$ is $F$-sunny, even simple).
Triangulate $X\x I$ and $Q\x I$ so that $F$ and the projections $P\:Q\x I\to Q$,
$p\:X\x I\to X$ are simplicial. Let $Y$ be the $(n-1)$-skeleton of $X$ and
$S=S(P\o F)$. By general position $\dim S\le n-1$, hence $S\su Y\x I$. Notice
that although the cylindric collapse $X\x I\col X\x 0\cup Y\x I$ is $F$-sunny
(even simple) and decreases the dimension, it, however, increases by 1 the
codimension of $S$ and the inductive step can not be applied. To overcome this,
we should arrange some device for collapsing away the top-dimensional simplexes
of $S$.

Let $A_i$ be a simplex of $S$. Since $p|_S$ is simplicial and non-degenerate,
there is a unique simplex $B$ of $Y$ such that $\Int A_i\su (\Int B)\x I$.
Since $A_i$ is locally of codimension $\ge 2$ in $X\x I$, $B$ faces some
$n$-simplex, say $C$, of $X$. Let $a$ be the barycenter of $A_i$. Suppose that
$A_i\not\su X\x 0$ and define some points near $a$. If $A_i\su X\x\Int I$,
choose $a_\uparrow$ directly above $a$ (with respect to $p$), $a_\downarrow$
directly below $a$ and $a_\to$ in $\Int(C\x I)$. If $A_i\su X\x 1$, let
$a_\uparrow=a$, choose $a_\downarrow$ directly below $a$ and $a_\to$ in
$\Int C$. Define a {\it blister} $J_i=a_\uparrow*a_\downarrow*a_\to*\partial
A_i$ over $A_i$ (here `$*$' means join), its {\it bad face}
$K_i=a_\uparrow*a_\downarrow*\partial A_i$ and its {\it good face}
$L_i=(a_\uparrow\cup a_\downarrow)*a_\to*\partial A_i$. For each simplex
$A_i$ of $S$ construct $J_i$ in so small neighborhood of $A_i$, that each
blister meets $Y\x I$ only in its bad face and
$J_i\cap J_j=\partial A_i\cap\partial A_j$ for each $i\neq j$. Let $J$, $K$ and
$L$ be the unions of all $J_i$, $K_i$ and $L_i$, respectively. Notice that
enough room to construct blisters is due to local codimension $\ge 2$ of $S$.

Since the blisters are small, we can collapse all the top-dimensional prisms,
but leaving the blisters. This gives an $F$-sunny collapse
$X\x I\col X\x 0\cup Y\x I\cup J$. Now each blister $J_i$ has the bad face
$K_i$ as a free face, which we may collapse it from. Therefore the blisters
$J_i$ can be collapsed onto their good faces in any order, particularly, in
the order the corresponding $A_i$'s $F$-overshadow each other. Although we
omitted simplexes, lying in $X\x 0$, the condition $F(X\x 0)\su Q\x 0$
implies that they $F$-overshadow nothing. Hence the obtained collapse $J\col L$
or, equivalently, $Y\x I\cup J\col (Y\x I\but K)\cup L=Z$, is $F$-sunny.
Thus there is an $F$-sunny collapse $X\x I\col Z\cup X\x 0$.

\bigskip
\def\epsfsize#1#2{\hsize}
\centerline{\epsffile{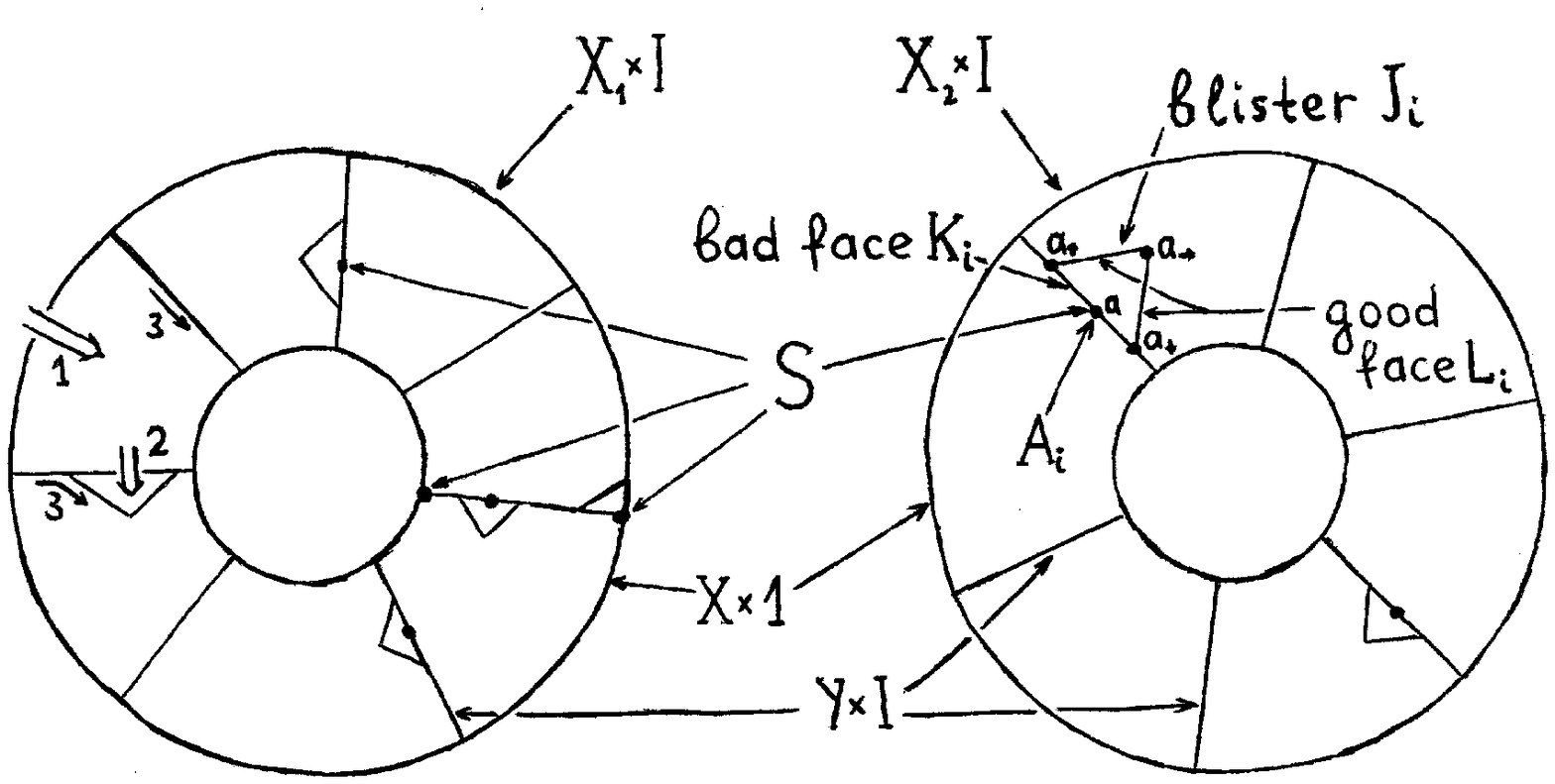}}
\medskip
\centerline{\bf Figure 2}
\bigskip

By mapping $a\mapsto a_\to$ for each $A_i$ and connecting linearly with
$\id|_{Z\cap Y\x I}$ we obtain a homeomorphism $\phi\:Y\x I\to Z$.
Notice that $Z$ meets its own $F$-shadow in $S'=\Cl{Z\cap S\but X\x 0}$ which
is, by the construction of blisters, locally of codimension $\ge 2$ in $Z$.
Also, $p|_{S'}=p|_{\phi^{-1}S'}$ is non-degenerate. Hence
$F\o\phi\:Y\x I\to Q\x I$ is in general position, and by the inductive
hypothesis $Y\x I$ collapses $(F\o\phi)$-sunny onto $Y\x 0$. Since $\phi$ is a
homeomorphism, $Z$ collapses $F$-sunny onto $Y\x 0$, hence we obtain an
$F$-sunny collapse $X\x I\col Z\cup X\x 0\col X\x 0$. \qed
\enddemo

\head 4. Lagging collapse \endhead

\definition{Definition} Let $F\:K\to Q\x I$ be a map. A collapse $V\col W$ in
$K$ is a {\it simple stable} $F$-sunny collapse if $V\but\Int W$ is not
$F$-overshadowed by $V$ (or, equivalently, $V\cap\sh_F V\su\Int W$; here
`$\Int$' denotes topological interior in $K$). Define {\it stable} $F$-sunny
collapse to be a sequence of simple stable collapses.
If $K=K_1\sqcup\dots\sqcup K_l$, let {\it $F$-link-shadow} of $V\su K$, denoted
$\lsh_F V$, be the union of $(\sh_F(V\cap K_i))\but K_i$ for all $i=1,\dots,l$
(speaking informally, link-shadow is the analogue of shadow for the case when
each $F(K_i)$ is `visible' from $F(K_j)$, $j\neq i$, but `transparent' from
$F(K_i)$). Define (simple) (stable) {\it $F$-link-sunny} collapse using
$F$-link-shadow instead of $F$-shadow. Evidently, any stable $F$-(link-)sunny
collapse is $F$-(link-)sunny, and any (stable) $F$-sunny collapse is (stable)
$F$-link-sunny, but not vice versa.
\enddefinition

\proclaim{Lemma 4.1} Let $X$ be a compact polyhedron, $Q$ a PL manifold and
$F\:X\x I\to Q\x I$ a PL link map such that $F(X\x 0)\su Q\x 0$. If there is an
$F$-link-sunny collapse $X\x I\col X\x 0$, then there is a stable
$F$-link-sunny collapse $X\x I\col X\x 0$.
\endproclaim

\proclaim{Sublemma 4.2}
(compare to \cite{12, Claim on p\. ~188}, \cite{13, Lem\. ~5.3})
Let $K=K_1\sqcup\dots\sqcup K_l$ be a finite simplicial complex, $Q$ a
combinatorial manifold, $Q\x I$ be triangulated so that the projection
$P\:Q\x I\to Q$ is simplicial, and $F\:K\to Q\x I$ be a simplicial link map.
Then there is a second derived subdivision $K''$ of $K$ such that
$\lsh_F V\su W$ implies $\lsh_F\N_{K''}(V)\su\Int\N_{K''}(W)$ for any
subcomplexes $V$ and $W$ of $K$.
\endproclaim

Here $\N_{K''}(L)$ denotes (for any subcomplex $L$ of $K$) the second derived
neighborhood of $L$ in $K$, which is the simplicial neighborhood of $L''$ in
$K''$.

\demo{Proof of Sublemma 4.2} Let $K'$ be the {\it barycentrically} derived
subdivision of $K$ and let us construct a derived subdivision $K''$ of
$K'$. Let $A_1,\dots,A_p$ be the simplexes of $K$, arranged in an order
of increasing dimension. Assuming that $A_i''$ is already defined for each
$i<j$, define $A_j''$ as follows. Let $a_j$ be the barycentre of $A_j$ and
$\Pi\:Q\x I\to I$ the horizontal projection. Define a map $f_{A_j}\:A_j\to\R^1$
by $\partial A_j\mapsto -1$, $a_j\mapsto\Pi\o F(a_j)+\frac{1}{100}$ and
extending linearly ($\frac{1}{100}$ can be replaced by any fixed positive
number). For any (already defined) derivation point $b$ of a simplex $B$ of
$(\partial A_j)'$, define a derivation point of $a_j*B$ to be
$(a_j*b)\cap f_{A_j}^{-1}(0)$.

Any point $x\in K$ is contained in the interior of a unique simplex $C_x$ of
$K$. If $c_x$ is the barycentre of $C_x$ and $x\neq c_x$, there is a unique
point $p_x\in\partial C_x$ such that $x\in p_x*c_x$. Notice that for any
subcomplex $L$ of $K$ the following criterion holds: $x\in\N_{K''}(L)\but L$
(respectively $x\in\Int\N_{K''}(L)\but L$) if and only if $x\notin L$,
$p_x\in\N_{K''}(L)$ (resp. $p_x\in\Int\N_{K''}(L)$) and $f_{C_x}(x)\le 0$
(resp. $f_{C_x}(x)<0$).

Suppose that a point $y\in K_i$ $F$-overshadows a point $z\in K_j$, $i\neq j$,
and prove the statement: `$y\in\N_{K''}(V)$ implies $z\in\Int\N_{K''}(W)$'
by induction on $\dim C_y$. If $\dim C_y=0$, $y$ is a vertex of $K$. Hence
$y\in\N_{K''}(V)$ implies $y\in V$, consequently $z\in W\su\Int\N_{K''}(W)$.
Let us prove the inductive step. Since $P$, $F$ are simplicial, $F(c_y)$
overshadows $F(c_z)$ and $F(p_y)$ either overshadows or equals $F(p_z)$; since
$i\neq j$ and $F$ is a link map, $F(p_y)\neq F(p_z)$. We can assume that
$y\notin V$, then $y\in\N_{K''}(V)$ implies $f_{C_y}(y)\le 0$. Since
$f_{C_z}(c_z)<f_{C_y}(c_y)$, we have that $f_{C_z}(z)<f_{C_y}(y)$, consequently
$f_{C_z}(z)<0$. Also $y\in\N_{K''}(V)$ implies $p_y\in\N_{K''}(V)$, hence
by the inductive hypothesys $p_z\in\Int\N_{K''}(W)$. Thus, if $z\notin W$,
we obtain $z\in\Int\N_{K''}(W)\but W$. \qed
\enddemo

\proclaim{Sublemma 4.3} If $V$ simplicially collapses onto $W$ in a finite
simplicial complex $K$, then $\N_{K''}(V)$ collapses onto $\N_{K''}(W)$.
\endproclaim

\demo{Proof} We will need the following simple observation. Suppose that a
simplex $A$ (strictly) faces a simplex $B$. Since ball collapses to its face,
$$\N_{B''}(A)\col\N_{\partial B''}(A)\cup\N_{B''}(\partial A).\tag{$*$}$$

Without loss of generality $\Cl{V\but W}$ is a single simplex, say $C$. Let $D$
be its free face, i.e\. $\Cl{\partial C\but W}$. Applying $(*)$ to $A=C$ and $B$
running over the faced by $C$ simplexes of $K$, in order of decreasing
dimension, we obtain $\N_{K''}(V)\col V\cup\N_{K''}(W\cup D)$. Applying
$(*)$ to $A=D$ and $B$ running over the faced by $D$ simplexes of $K$ except
for $C$, in order of decreasing dimension, we obtain
$V\cup\N_{K''}(W\cup D)\col V\cup\N_{K''}(W)$. Finally, since ball
collapses to its face, $V\cup\N_{K''}(W)\col\N_{K''}(W)$. \qed
\enddemo

\def\epsfsize#1#2{\hsize}
\centerline{\epsffile{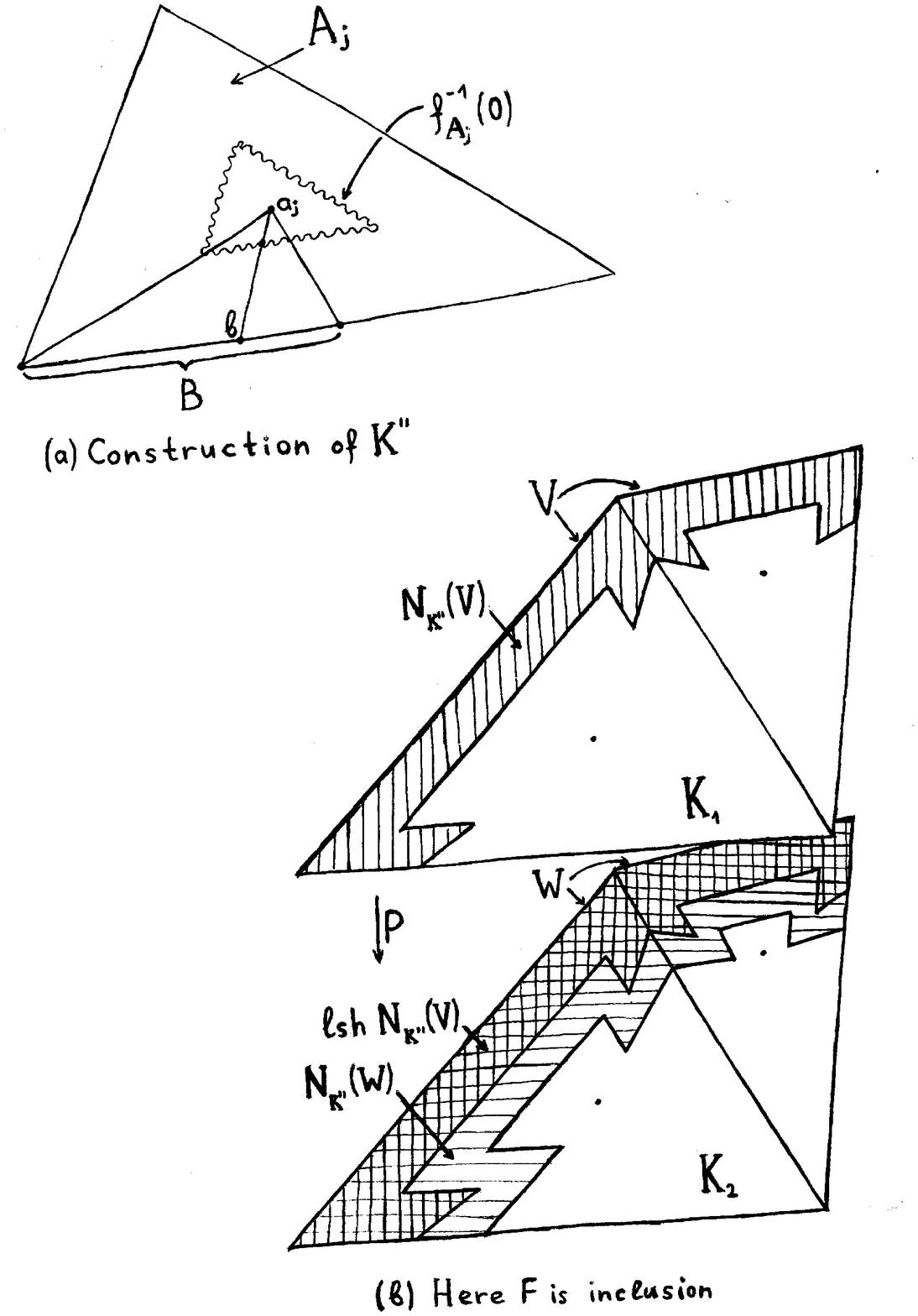}}
\medskip
\centerline{\bf Figure 3}
\bigskip

\demo{Proof of Lemma 4.1} We are given collapses
$X\x I=K_0\col\dots\col K_m=X\x 0$ such that $K_i\cap\lsh_F K_i\su K_{i+1}$.
By induction on $i$ it follows that $\lsh_F K_i\su K_{i+1}$. Triangulate
$X\x I$ so that $K_1,\dots,K_m$ are its subcomplexes, and both $P$ and $F$ are
simplicial in some triangulation of $Q\x I$. Let $U_i=\N_{(X\x I)''}(K_i)$,
$i=0,\dots,m$, where the subdivision $(X\x I)''$ is yielded by 4.2. Then
$\lsh_F U_i\su\Int U_{i+1}$ (hereafter `Int' means topological interior in
$X\x I$). Furthermore, $F(X\x 0)\su Q\x 0$ implies that
$\N_{(X\x I)''}(X_j\x 0)$ can $F$-overshadow only the points of $X_j\x I$.
Hence the $F$-link-shadow of $U_m$ is empty. By 4.3, $U_i\col U_{i+1}$, and
since $U_m$ is a regular neighborhood of $X\x 0$ (or by an argument analogous
to the proof of 4.3), $U_m\col X\x 0$. Thus we obtain a stable $F$-link-sunny
collapse $X\x I=U_0\col\dots\col U_m\col X\x 0$. \qed
\enddemo

\head 5. Proofs of 1.1 and 1.3 \endhead

\demo{Proof of Theorem 1.1} Let $F\:X\x I\to Q\x I$ be the given link
concordance between $f_0$ and $f_1$. Without loss of generality $F$ is PL and, by
2.1, in general position. Then 3.1 and 4.1 yield a stable $F$-link-sunny
collapse $X\x I\col X\x 0$. This collapse gives a homotopy $h_t\:X\x I\to X\x I$
(obtained by retracting, for each elementary collapse $V\col W$, of the ball
$\Cl{V\but W}$ onto its face $\Cl{V\but W}\cap W$). Conversely, for each $t\in I$
the map $h_t$ corresponds to some simple stable $F$-link-sunny collapse
$V\col W$. Then for any $Y\su X\x I$ we have $h_t(Y)\su V$ and
$h_t(Y)\cap\Int W=Y\cap\Int W$. Consequently,
$h_t(X\x 1)\su V\but ((\Int W)\but X\x 1)$. On the other hand, since $V\col W$ is
a simple stable $F$-link-sunny collapse and $F(X\x 1)\su Q\x 1$, we have that
$\lsh_F(V)\su(\Int W)\but X\x 1$. Thus $h_t(X\x 1)$ does not $F$-link-overshadow
itself.

Define $H\:X\x I\to X\x I$ by $H(x,t)=h_{2t}(x\x 1)$ for $t\in [0,\frac{1}2]$
and $H(x,t)=p\o h_{2-2t}(x\x 1)$ for $t\in [\frac{1}2,1]$ (where
$p\:X\x I\to X$ is the projection). Then $H$ fixes $X\x\{0,1\}$, while
$H(X\x t)$ does not $F$-link-overshadow itself. In other words,
$F\o H(X_i\x t)$ and $F\o H(X_j\x t)$ do not overshadow each other whenever
$i\neq j$. Since $F$ is a link map, they do not intersect as well. Hence
a map $\Phi\:X\x I\to Q\x I$, defined by $\Phi(x,t)=(P\o F\o H(x,t),t)$, is a
link homotopy between $f_0$ and $f_1$. \qed
\enddemo

\demo{Proof of Theorem 1.3} Using the above notation, let us show that if $F$
is a doodle concordance, then $\Phi$ is actually a doodle homotopy.
By the above, $F\o H(X_i\x t)$ and $F\o H(X_{j_1}\x t)\cap F\o H(X_{j_2}\x t)$
do not overshadow each other for any distinct $i,j_1,j_2$. Since $F$ is a
doodle, they do not meet as well. Therefore $\Phi(X_i\x t)$ does not meet
$\Phi(X_{j_1}\x t)\cap\Phi(X_{j_2}\x t)$. It follows that $\Phi$ is a doodle
homotopy. \qed
\enddemo

\head 6. Proof of 1.4 \endhead

\definition{Definition} Let $K$ be a polyhedron and $F\:K\to Q\x I$ a map.
Define {\it $F$-$\eps$-shadow} of $V\su K$ to be
$\sh_F(V)\but\N_\eps(F^{-1}F(V))$ and define {\it (simple) (stable)
$F$-$\eps$-sunny} collapse using $F$-$\eps$-shadow instead of $F$-shadow.
The proof of the following lemma is analogous to that of Lemma 4.1.
\enddefinition

\proclaim{Lemma 6.1} For any $\eps>0$ there is $\dta>0$ such that the following
holds. Let $X$ be a compact polyhedron, $Q$ a PL manifold and
$F\:X\x I\to Q\x I$ a PL map such that $F(X\x 0)\su Q\x 0$. If there is an
$F$-$\dta$-sunny collapse $X\x I\col X\x 0$, then there is a stable
$F$-$\eps$-sunny collapse $X\x I\col X\x 0$. \qed
\endproclaim

The below statements follow immediately from the proofs of 3.1 and 4.1.

\proclaim{Addendum 6.2 to Lemma 3.1} For each $\eps>0$ we can choose the
$F$-sunny collapse $X\x I\col X\x 0$ so that the trace of $x\x I$ lies in
$\N_\eps(x)\x I$ for any $x\in X$. \qed
\endproclaim

\proclaim{Addendum 6.3 to Lemma 6.1} The stable $F$-$\eps$-sunny collapse
$X\x I\col X\x 0$ can be chosen so that the trace of any point $p\in X\x I$
under this collapse lies in the $\eps$-neighborhood of that under the given
collapse. \qed
\endproclaim

\demo{Proof of Theorem 1.4} Without loss of generality the given map
$F\:X\x I\to Q\x I$ is a PL general position map. Apply (3.1+6.2) and
(6.1+6.3) to obtain a stable $F$-$\frac{\eps}2$-sunny collapse
$X\x I\col X\x 0$ such that the trace of $x\x I$ lies in $\N_\heps(x)\x I$ for
any $x\in X$. Analogously to the proof of 1.1 we obtain a map
$H\:X\x I\to X\x I$ fixing $X\x \{0,1\}$ and such that

1) $H(x\x I)\su\N_\heps(x)\x I$ for each $x\in X$, and

2) $X\x t$ does not $(F\o H)$-$\frac{\eps}2$-overshadow itself for each
$t\in I$.

\noindent
Let us prove that $\Phi\:X\x I\to Q\x I$, defined by
$\Phi(x,t)=(P\o F\o H(x,t),t)$, is the required map. Indeed, suppose that
$\Phi(x_1\x I)\cap\dots\cap\Phi(x_l\x I)\neq\emptyset$ for some
$x_1,\dots,x_l\in X$. Since $\Phi$ is level-preserving,
$\Phi(x_1\x t)=\dots=\Phi(x_l\x t)$ for some $t\in I$. Then
$F\o H(x_1\x t),\dots,F\o H(x_l\x t)$ lie in the same vertical line in $Q\x I$.
By (2) $F\o H(\N_\heps(x_1)\x t)\cap\dots\cap F\o H(\N_\heps(x_l)\x t)\neq
\emptyset$. By (1) $F\o H(x_i\x I)\su F(\N_\heps(x_i)\x I)$ for $i=1,\dots,l$.
Hence $F(\N_\eps(x_1)\x I)\cap\dots\cap F(\N_\eps(x_l)\x I)\neq\emptyset$. \qed
\enddemo

\head Acknowledgements \endhead

I wish to thank to A. Skopenkov for suggesting the problem, encouragment
and guidance in preparation of this paper, to U. Koschorke for sending his
works on the subject, to P. M. Akhmetyev, V. M. Nezhinsky, R. Sadykov,
K. Salikhov and E. V. \v{S}\v{c}epin for useful discussions and important
remarks.

\head Note added in 2018 \endhead

This preprint dates from 1998, when it was privately circulated and publicly presented
(Pontryagin's 90th Anniverary Conference, Moscow, September 1998;
Rokhlin Topology Seminar of V. M. Nezhinsky, St.-Petersburg, March 1998;
Geometric Topology Seminar of E. V. Shchepin, Moscow, Fall 1997).
Being a second year undergraduate student at the time, I followed the advice by
A.~ Skopenkov to submit the preprint to the journal Topology --- where it was 
received in April 1998 and rejected in September 2000.

The results of this preprint were announced in Uspekhi Mat.\ Nauk 55:3 (2000), 
183--184 (English transl.: Russ.\ Math.\ Surv.\ 55 (2000), 589--590).
A more elaborate verion of the main construction (the ``lagging collapse''), which 
reproves the classical Concordance Implies Isotopy Theorem and proves its controlled 
version, appeared in Topol.\ Appl.\ 120 (2002), 105--156; arXiv:\,math.GT/0101047 (\S 5).
Given that, and the fact that Theorem 1.1 is far from being optimal (as compared to 
the results of X.-S. Lin and P. Teichner) I was no longer keen to publicize 
the present preprint in its extant form.
However, the referee of another paper, which uses Theorem 1.3, leaves me no choice.


\Refs \refstyle{C}

\ref \no 1 \by R. H. Bing
\paper Vertical general position
\inbook Geometric topology (proceedings of conference, Utah),
LNiM {\bf 438} \publ Springer \yr 1975 \pages 16--42
\endref

\ref \no 2 \by R. D. Edwards
\paper The equivalence of close piecewise-linear embeddings
\jour Gen. Topol. and Appl. \vol 5 \yr 1975 \pages 147--180
\endref

\ref \no 3 \by R. Fenn, P. Taylor
\paper Introducing doodles
\inbook Topology of low-dimensional manifolds (proceedings of II Sussex
conference), LNiM {\bf 722} \ed R. Fenn \publ Springer \publaddr Berlin
\yr 1979 \pages 37--43
\endref

\ref \no 4 \by R. Fenn, D. Rolfsen
\paper Spheres may link homotopically in 4-space
\jour J. Lond. Math. Soc. \vol 34 \yr 1986 \pages 177--184
\endref

\ref \no 5 \by C. H. Giffen
\paper Link concordance implies link homotopy
\jour Math. Scand. \vol 45 \yr 1979 \pages 243--254
\endref

\ref \no 6 \by D. L. Goldsmith
\paper Concordance implies homotopy for classical links in $M^3$
\jour Comm. Math. Helv. \vol 54 \yr 1979 \pages 347--355
\endref

\ref \no 7 \by N. Habegger
\paper On linking coefficients
\jour Proc. AMS \vol 96 \yr 1986 \pages 353--359
\endref

\ref \no 8 \bysame
\paper Applications of Morse theory to link theory \inbook Knots 90
\publ de Gruyter \publaddr Berlin, N.Y. \yr 1992 \pages 389--394
\endref

\ref \no 9 \by N. Habegger, U. Kaiser
\paper Link homotopy in the 2-metastable range
\jour Topology \vol 37 \yr 1998 \pages 75--94
\endref

\ref \no 10 \by A. Haefliger
\paper Differentiable embeddings of $S^n$ in $S^{n+q}$ for $q>2$
\jour Ann. Math. \vol 83 \yr 1966 \pages 402--436
\endref

\ref \no 11 \by L. S. Harris
\paper Intersections and embeddings of polyhedra
\jour Topology \vol 8 \yr 1969 \pages 1--26
\endref

\ref \no 12 \by J. F. P. Hudson
\book Piecewise-linear topology
\publ Benjamin W. A. \publaddr N. Y., Amsterdam \yr 1969
\endref

\ref \no 13 \bysame
\paper Concordance, isotopy and diffeotopy
\jour Ann. Math. \vol 91 \yr 1970 \pages 425--448
\endref

\ref \no 14 \by J. F. P. Hudson, W. B. R. Lickorish
\paper Extending piecewise-linear concordances
\jour Quart. J. Math. \vol 22 \yr 1971 \pages 1--12
\endref

\ref \no 15 \by M. Khovanov
\paper Doodle groups
\jour Trans. AMS \vol 349 \yr 1997 \pages 2297--2315
\endref

\ref \no 16 \by P. A. Kirk
\paper Link homotopy with one codimension 2 component
\jour Trans. AMS \vol 319 \yr 1990 \pages 663--688
\endref

\ref \no 17 \by U. Koschorke
\paper Link maps and the geometry of their invariants
\jour Manuscripta Math. \vol 61 \yr 1988 \pages 383--415
\endref

\ref \no 18 \bysame
\paper On link maps and their homotopy classification
\jour Math. Ann. \vol 286 \yr 1990 \pages 753--782
\endref

\ref \no 19 \bysame
\paper Link homotopy with many components
\jour Topology \vol 30 \yr 1991 \pages 267--281
\endref

\ref \no 20 \bysame
\paper Homotopy, concordance and bordism of link maps
\inbook Global analysis in modern mathematics \publ Publish or Perish, Inc.
\publaddr Houston/Texas \yr 1993 \pages 283--299
\endref

\ref \no 21 \bysame
\paper A generalization of Milnor's $\mu$-invariants to higher dimensional
link maps \jour Topology \vol 36 \yr 1997 \pages 301--324
\endref

\ref \no 22 \by U. Koschorke, D. Rolfsen
\paper Higher dimensional link operations and stable homotopy
\jour Pacific J. Math. \vol 139 \yr 1989 \pages 87--106
\endref

\ref \no 23 \by W. B. R. Lickorish
\paper The piecewise-linear unknotting of cones
\jour Topology \vol 4 \yr 1965 \pages 67--91
\endref

\ref \no 24 \by X. S. Lin
\paper On equivalence relations of links in 3-manifolds
\publaddr UCSD \miscnote preprint \yr 1985
\endref

\ref \no 25 \by W. S. Massey, D. Rolfsen
\paper Homotopy clasification of higher dimensional links
\jour Indiana Univ. Math. J. \vol 34 \yr 1986 \pages 375--391
\endref

\ref \no 26 \by J. W. Milnor
\paper Link groups
\jour Ann. Math. \vol 59 \yr 1954 \pages 177--195
\endref

\ref \no 27 \by K. Millett
\paper Piecewise-linear concordances and isotopies
\jour Mem. AMS \vol 153 \yr 1974 \pages 1--73
\endref

\ref \no 28 \by V. M. Nezhinskij
\paper Groups of classes of pseudohomotopic singular links. I
\jour Zapiski nauchnykh seminarov LOMI \vol 168 \yr 1988 \pages 114--124
\lang in Russian
\endref

\ref \no 29 \bysame
\paper Singular links of type $(p,2k+1)$ in the $(4k+2)$-dimensional sphere
\jour Algebra i Analis \vol 5 \yr 1993 \pages 170--190 \lang in Russian
\moreref \paper II
\jour Zapiski nauchnykh seminarov LOMI \vol 231 \yr 1995 \pages 191--196
\lang in Russian
\endref

\ref \no 30 \by C. P. Rourke
\paper Embedded handle theory, concordance and isotopy
\inbook Topology of Manifolds \eds J. C. Cantrell, C. H. Edwards
\publ Markham \publaddr Chicago, Ill. \yr 1970 \pages 431--438
\endref

\ref \no 31 \by C. P. Rourke, B. J. Sanderson
\book Introduction into piecewise-linear topology
\publ Springer-Verlag \publaddr Berlin, N.Y. \yr 1972
\endref

\ref \no 32 \by R. F. Sayakhova
\paper Singular links of type $(1,1,m)$ in $S^3$, $m>1$
\miscnote preprint \lang in Russian
\endref

\ref \no 33 \by G. P. Scott
\paper Homotopy links
\jour Abh. Math. Sem. Univ. Hamburg \vol 32 \yr 1968 \pages 186--190
\endref

\ref \no 34 \by J. Segal,  A. B. Skopenkov, S. Spie\D z
\paper Embeddings of polyhedra in $\R^m$ and the deleted product obstruction
\jour Topol. Appl. \vol 85 \yr 1997 \pages 1--10
\endref

\ref \no 35 \by A. B. Skopenkov
\paper On the generalized Massey--Rolfsen invariant for link maps
\miscnote preprint
\endref

\ref \no 36 \by P. Teichner
\paper Gropes, Alexander duality and link homotopy
\jour Geometry and Topology \vol 1 \yr 1997 \pages 51--69
\endref

\ref \no 37 \by A. G. Tristram
\paper Some cobordism invariants for links
\jour Proc. Camb. Phil. Soc. \vol 66 \yr 1969 \pages 251--264
\endref

\ref \no 38 \by E. C. Zeeman
\paper Unknotting combinatorial balls
\jour Ann. Math. \vol 78 \yr 1963 \pages 501--526
\endref

\ref \no 39 \bysame
\paper Relative simplicial approximation
\jour Proc. Camb. Phil. Soc. \vol 60 \yr 1964 \pages 39--41
\endref

\endRefs
\enddocument
\end